\documentclass[11pt]{amsart}

\usepackage[margin=1.25in]{geometry}

\usepackage{amssymb}
\usepackage{amsmath}
\usepackage{mathrsfs}
\usepackage{url}
\usepackage{verbatim}
\input xy
\xyoption{all}
\usepackage{graphicx}  
\usepackage{tikz}
\usepackage{tkz-graph}
\usepackage{multicol}
\usepackage{enumerate}
\usepackage{mathtools}
\usepackage{url}
\usepackage{hyperref}

\newtheorem{thm}{Theorem}

\newtheorem{cor}[thm]{Corollary}
\newtheorem{prop}[thm]{Proposition}

\newtheorem{Definition}[thm]{Definition}
\newenvironment{definition}
  {\begin{Definition}\rm}{\end{Definition}}

\newtheorem{Example}[thm]{Example}

\newtheorem{Remark}[thm]{Remark}

\usepackage[enableskew]{youngtab}

\usepackage{mathdots}

\newcommand\subdot{\mathrel{\ooalign{$\subseteq$\cr
  \hidewidth\hbox{\scalebox{1.5}{$\cdot\mkern1mu$}}\cr}}}
  
\allowdisplaybreaks

\title{A note on statistical averages for oscillating tableaux}
\author{Sam Hopkins}
\email{shopkins@mit.edu}
\author{Ingrid Zhang}
\email{ingridz@mit.edu}
\address{Massachusetts Institute of Technology, Cambridge MA, 02139}

\begin{document}

\begin{abstract}
We define a statistic called the weight of oscillating tableaux. Oscillating tableaux, a generalization of standard Young tableaux, are certain walks in Young's lattice of partitions. The weight of an oscillating tableau is the sum of the sizes of all the partitions that it visits. We show that the average weight of all oscillating tableaux of shape $\lambda$ and length $|\lambda|+2n$, where $|\lambda|$ is the size of $\lambda$, has a surprisingly simple formula: it is a quadratic polynomial in $|\lambda|$ and $n$. Our proof via the theory of differential posets is largely computational. We suggest how the homomesy paradigm of Propp and Roby may lead to a more conceptual proof of this result and reveal a hidden symmetry in the set of perfect matchings.
\end{abstract}

\maketitle

In this note, we follow the standard notation for partitions as laid out in~\cite[\S7.2]{stanley2}. Recall that Young's lattice is the poset of all partitions partially ordered by inclusion of Young diagrams. We write $\mu \subdot \lambda$ to mean that $\lambda$ covers $\mu$ in Young's lattice; that is, $\mu \subdot \lambda$ means that the Young diagram of $\lambda$ is obtained from that of $\mu$ by the addition of a box. A \emph{walk in the (Hasse diagram of) Young's lattice} is a sequence of partitions $(\lambda^{0},\lambda^{1},\ldots,\lambda^{l})$ such that for all $1 \leq i \leq l$ we have $\lambda^{i-1} \subdot \lambda^{i}$ or~$\lambda^{i} \subdot \lambda^{i-1}$. An \emph{oscillating tableau of shape $\lambda$ and length $l$} is a walk $T = (\lambda^{0},\lambda^{1},\ldots,\lambda^{l})$ in Young's lattice such that $\lambda^{0} = \emptyset$ and $\lambda^{l} = \lambda$. We use $\mathcal{OT}(\lambda,l)$ to denote the set of oscillating tableaux of shape $\lambda$ and length~$l$.

Oscillating tableaux are a generalization of standard Young tableaux (SYT): it is not hard to see that $\#\mathcal{OT}(\lambda,|\lambda|) = f^{\lambda}$, the number of SYT of shape $\lambda$. Here $|\lambda|$ is the size of $\lambda$. The numbers $f^{\lambda}$ are easily computed thanks to the famous hook-length formula~\cite[Corollary 7.21.6]{stanley2}. In fact, all oscillating tableaux enjoy simple enumerative formulas.

\begin{thm} \label{thm:count}
Let $|\lambda|=k$. Then for all $n \in \mathbb{N}$,
\[\#\mathcal{OT}(\lambda,k+2n) = \binom{2n+k}{k}(2n-1)!! f^{\lambda}.\]
On the other hand, if $l \neq k + 2n$ for some $n \in \mathbb{N}$ then $\#\mathcal{OT}(\lambda,l) = 0$.
\end{thm}

Above, $(2n-1)!! := 1\cdot3\cdot5\cdots(2n-1)$, which is $1$ if $n=0$. The last sentence in Theorem~\ref{thm:count} is clear because for any walk $(\lambda^{0},\ldots,\lambda^{l})$ in Young's lattice we have $|\lambda^{l}| - |\lambda^{0}| = l \mod 2$ and $|\lambda^{l}| \leq |\lambda^{0}| + l$. There are a variety of approaches to proving the nontrivial claim of this theorem. A bijective proof is given in~\cite[Lemma 8.7]{sundaram}. Roby~\cite[\S4.2]{roby} explains how the theory of growth diagrams, initiated by Fomin~\cite{fomin}, also yields a bijective proof. Another elegant approach to this formula is via the theory of differential posets developed by Stanley~\cite{stanley1}~\cite[Theorem 8.7]{stanley3}; this is the approach we will follow in investigating certain statistical averages for oscillating tableaux.

Oscillating tableaux have, to some extent, been studied in the context of combinatorial statistics; for example, Chen~et.~al.~\cite{chen} used tableaux to prove the symmetry of maximal crossing number and maximal nesting number in set partitions. Here we are interested in a statistic on tableaux, which we term weight, that, to our knowledge, has not been explored in previous research. It turns out that there is a surprisingly simple formula for the average weight across all oscillating tableaux of given shape and length. For~$T \in \mathcal{OT}(\lambda,l)$, define the \emph{weight} of $T$, denoted $\mathrm{wt}(T)$, to be the sum of the sizes of all partitions that $T$ visits; that is, if $T = (\lambda^{0},\ldots,\lambda^{l})$ then we set~$\mathrm{wt}(T) := \sum_{i=0}^{l} |\lambda^{i}|$. 

\begin{thm} \label{thm:main}
Let $|\lambda| = k$. Then for all $n \in \mathbb{N}$,
\[  \frac{1}{\# \mathcal{OT}(\lambda,k+2n)} \sum_{T \in \mathcal{OT}(\lambda,k+2n)} \mathrm{wt}(T) = \frac{1}{6}(4n^2 + 3k^2 + 8 kn + 2n +3k).\]
\end{thm}

It also makes sense to divide the average in Theorem~\ref{thm:main} by $(2n+k+1)$ in order to compute the average of $|\lambda^{i}|$ among all $i = 0,\ldots,k+2n$ and all $(\lambda^{0},\ldots,\lambda^{k+2n}) \in \mathcal{OT}(\lambda,k+2n)$. The following corollary asserts a clean formula for this average size.

\begin{cor} \label{cor:main}
Let $|\lambda| = k$. Then for all $n \in \mathbb{N}$,
\[  \frac{1}{\# \mathcal{OT}(\lambda,k+2n)\cdot(2n+k+1)} \sum_{T \in \mathcal{OT}(\lambda,k+2n)} \mathrm{wt}(T) = \frac{n}{3} + \frac{k}{2}.\]
\end{cor}

The fact that the average in Theorem~\ref{thm:main} depends only on $k$ and $n$, and not $\lambda$, is understandable in light of Theorem~\ref{thm:count}, which tells us that the dependence of~$\#\mathcal{OT}(\lambda,k+2n)$ on~$\lambda$ is slight. Still, it is unclear a priori why this average should be a polynomial in $k$ and $n$. Indeed, it is not clear that there should be an absolute constant $\kappa \in \{1,2,\ldots\}$ such that $\kappa$ times this average is integral for all $\lambda$ and all $n$. Note that because
\[ \frac{1}{2}(4n^2 + 3k^2 + 8 kn + 2n +3k) =(n+3k/2)(2n+k+1)  = (n+(k+1)/2)(2n+3k),\]
which is evidently integral for all $n,k \in \mathbb{N}$, we may take $\kappa = 3$. It is especially unclear why such a small choice of $\kappa$ should work. The proof of Theorem~\ref{thm:main} we give in~\S\ref{sec:proof}, essentially a mysterious computation using differential poset operators, is unilluminating in this regard. In~\S\ref{sec:homomesy} we suggest how the recent homomesy paradigm of Propp and Roby~\cite{propp} might lead to a more explanatory proof of Theorem~\ref{thm:main}. Homomesy says, roughly, that ``small denominators in statistical averages of combinatorial objects should be explained by group actions.'' The possibility of such a proof points towards a hidden $\mathfrak{S}_3$-symmetry in the set of perfect matchings, as well as potentially the first instance of homomesy involving a non-cyclic group.

One might hope to extend Theorem~\ref{thm:main} to skew oscillating tableaux. A \emph{skew oscillating tableau of shape $\lambda/\mu$ and length $l$} is a walk $T = (\lambda^{0},\lambda^{1},\ldots,\lambda^{l})$ in Young's lattice such that $\lambda^{0} = \mu$ and $\lambda^{l} = \lambda$. See~\cite[Corollary 4.2.11]{roby} for a formula enumerating such tableaux. However, computation carried out with Sage mathematical software~\cite{sage-combinat} suggests that there is not as simple a formula for the average weight of such tableaux for arbitrary~$\mu$; in particular, it appears that the denominator of the average is unbounded in this case.

\noindent {\bf Acknowledgements}: This research was carried out at MIT as part of the RSI summer mathematics research program for high-school students. The first author was the mentor of the second author. We thank Richard Stanley for some helpful questions and comments. We also thank Tanya Khovanova for proofreading and comments about the exposition.

\section{Proof of Theorem~\ref{thm:main}} \label{sec:proof}

To prove Theorem~\ref{thm:main}, we apply the theory of differential posets~\cite{stanley1}~\cite[Chapter 8]{stanley3}. Recall the general set up: we let $V$ be the vector space of linear combinations of partitions and define two linear operators $U,D\colon V \to V$ by~$U(\mu) := \sum_{\mu \subdot \nu} \nu$ and~$D(\mu) := \sum_{\nu \subdot \mu} \nu$. These operators (which multiply right-to-left) satisfy the fundamental identity $DU - UD = I$. Following Stanley~\cite[Chapter 8]{stanley3}, define $b_{ij}(l) \in \mathbb{Z}$ by $\left(U + D\right)^{l} = \sum_{i,j} {b}_{ij}(l) U^iD^j$. Lemma 8.6 of~\cite{stanley3} shows that $b_{i0}(l) = \binom{l}{i} (l-i-1)!!$ when $l-i$ is even. We basically mimic the proof of Lemma 8.6, but we expand a slightly more complicated expression that keeps track of tableaux weights. So let us similarly define $p_{ij}(l) \in \mathbb{Z}[x,x^{-1},y,y^{-1}]$ by
\[ \left(y^{x \, d/dx}xU + y^{x \, d/dx}x^{-1}D\right)^{l} = \sum_{i,j} p_{ij}(l) U^iD^j.\]
where $y^{x \, d/dx}\colon \mathbb{Z}[x,x^{-1}] \to \mathbb{Z}[x,x^{-1},y,y^{-1}]$ is the $\mathbb{Z}$-linear map of rings $y^{x \, d/dx}(x^n) := x^ny^n$. It is clear that the $p_{ij}(l)$ exist and are well defined. It is also clear that $p_{ij}(l) = x^{i-j}q_{ij}(l)$ for some unique $q_{ij}(l) \in \mathbb{Z}[y,y^{-1}]$. Set $\mathcal{T} :=   \mathcal{OT}(\lambda,k+2n)$. Then $\#\mathcal{T} = [\lambda] \left(U + D\right)^{k+2n} \cdot \, \emptyset.$ Similarly, we have
\[ \sum_{T \in \mathcal{T}} \mathrm{wt}(T) = [\lambda] \frac{d}{dy}\left(y^{x \, d/dx}xU + y^{x \, d/dx}x^{-1}D\right)^{k+2n}\mid_{x=1,y=1} \cdot \, \emptyset. \]
In this expression, the power of $x$ keeps track of the number of $U$s minus the number of~$D$s reading right-to-left in each term in the expansion of $(U+D)^{k+2n}$ and the power of $y$ keeps track of the running sum of the powers of $x$s. In other words, $x$ keeps track of the size of the partition we are at and $y$ keeps track of the weight of the tableau built up so far. Thus, defining $c_{ij}(l) \in \mathbb{Z}$ by~$c_{ij}(l) := \frac{d}{dy} (q_{ij}(l))\mid_{y=1}$, we have $\# \mathcal{T} = b_{k,0}(k+2n)f^{\lambda}$ and~$\sum_{T \in \mathcal{T}} \mathrm{wt}(T) = c_{k,0}(k+2n) f^{\lambda}$. The theorem is therefore equivalent to
\begin{equation}
\frac{c_{k0}(k+2n)}{b_{k0}(k+2n)} =  \frac{1}{3}\left(n+\frac{k+1}{2}\right)(2n+3k).  \label{eqn:key}
\end{equation}
Note that we have $DU^{i} = U^{i}D + iU^{i-1}$ for all $i \geq 0$, which follows inductively from the identity $DU - UD = I$. Using this fact we get
\begin{align*}
\sum_{i,j} x^{i-j}q_{ij}(l+1) U^iD^j &=  \left(y^{x \, d/dx}xU + y^{x \, d/dx}x^{-1}D\right)^{l}  \\
&= \left(y^{x \, d/dx}xU + y^{x \, d/dx}x^{-1}D\right) \sum_{i,j} x^{i-j} q_{ij}(l) U^iD^j \\
&= \sum_{i,j} x^{i-j+1}y^{i-j+1}q_{ij}(l) U^{i+1}D^j +  x^{i-j-1}y^{i-j-1}q_{ij}(l) DU^{i}D^j \\
&= \sum_{i,j} x^{i-j} y^{i-j}(q_{i-1,j}(l) + q_{i,j-1}(l-1) + (i+1)q_{i+1,j}(l-1)) U^{i}D^{j}
\end{align*}
so by equating the coefficients of $U^iD^j$ on both sides we get that the polynomials $q_{ij}(l)$ satisfy the recurrence
\[ q_{ij}(l+1) = y^{i-j}(q_{i-1,j}(l) + q_{i,j-1}(l) + (i+1)q_{i+1,j}(l)).\]
In particular, $q_{i0}(l+1) = y^{i}(q_{i-1,0}(l) + (i+1)q_{i+1,j}(l))$. Observe that~$q_{ij}(l)\mid_{y=1} = b_{ij}(l)$. So we get the following recurrence in the $c$ numbers:
\[ c_{i0}(l+1) = c_{i-1,0}(l) + (i+1)c_{i+1,0}(l)+ib_{i-1,0}(l) + i(i+1)b_{i+1,0}(l).\]
Assume that~(\ref{eqn:key}) holds for smaller values of $k+2n$. Then
\begin{align*}
c_{k0}(k+2n) = \, &c_{k-1,0}(k-1+2n) +kb_{k-1,0}(k-1+2n) \\
&+ (k+1)c_{k+1,0}(k+1+2(n-1)) + k(k+1)b_{k+1,0}(k-1+2n)\\
=&\left(\frac{(n+\frac{k}{2})(2n+3(k-1))}{3}+k\right)b_{k-1,0}(k-1+2n) \\
&+  \left(\frac{(n-1+\frac{k+2}{2})(2(n-1)+3(k+1))}{3}  + k\right)(k+1)b_{k+1,0}(k-1+2n) \\
=& \left( \frac{(n+\frac{k+1}{2})(2n+3k)}{3} -\frac{4}{3}n\right)b_{k-1,0}(k-1+2n) \\
&+ \left(\frac{(n+\frac{k+1}{2})(2n+3k)}{3} +\frac{2}{3}k\right)(k+1)b_{k+1,0}(k-1+2n) \\
=& \left(\frac{(n+\frac{k+1}{2})(2n+3k)}{3}\right) (b_{k-1,0}(k-1+2n) +(k+1)b_{k+1,0}(k-1+2n) ) \\
&+ \frac{2}{3} \Big(k(k+1)b_{k+1,0}(k-1+2n) -2nb_{k-1,0}(k-1+2n)\Big) \\
=& \left(\frac{(n+\frac{k+1}{2})(2n+3k)}{3}\right) b_{k,0}(k+2n) \\
&+ \frac{2}{3} \left( k(k+1)\binom{k-1+2n}{k+1}(2n-3)!! - 2n\binom{k-1+2n}{k-1}(2n-1)!!\right) \\
=& \left(\frac{(n+\frac{k+1}{2})(2n+3k)}{3} \right) b_{k,0}(k+2n).
\end{align*}
So~(\ref{eqn:key}) holds for $k$ and $n$. Also $c_{00}(0) = 0$, which agrees with~(\ref{eqn:key}). Thus we are done by induction.

\section{A better proof through homomesy?} \label{sec:homomesy}

In the study of statistics on combinatorial objects, Propp and Roby~\cite{propp} recently introduced the paradigm of homomesy.

\begin{definition}
Let $\mathbf{k}$ be a field of characteristic zero. We say that the statistic $f:\mathcal{S} \to \mathbf{k}$ on a set of combinatorial objects $\mathcal{S}$ is \emph{homomesic} with respect to the action of a group $G$ on~$\mathcal{S}$ if there is some $c \in \mathbf{k}$ such that~$\frac{1}{\#\mathcal{O}} \sum_{s \in \mathcal{O}} f(s) = c$ for each $G$-orbit $\mathcal{O}$. In other words, we say $f$ is homomesic with respect to $G$ if the average of $f$ is the same for each $G$-orbit. In this case we also say the triple $(\mathcal{S},G,f)$ exhibits homomesy.
\end{definition}

Propp and Roby point out how, like the cyclic-sieving phenomenon,  the homomesy phenomenon is widespread in algebraic and enumerative combinatorics. Past research on homomesy has in general restricted to the case $G = \langle \varphi \rangle$ where $\varphi\colon \mathcal{S} \to \mathcal{S}$ is some combinatorially-meaningful procedure such as rowmotion in order ideals of posets~\cite{propp}~\cite{einstein} or promotion of Young tableaux~\cite{bloom}. The goal in such settings is to find natural statistics on $\mathcal{S}$ that are homomesic with respect to $\varphi$. However, we could take a different perspective whereby we have a certain statistic $f$ we want to understand better and search for a group $G$ so that~$(\mathcal{S},f,G)$ exhibits homomesy. Of course, this is always possible by taking $G$ to be the group of all invertible maps $\mathcal{S} \to \mathcal{S}$; but this is a trivial case in which there is only one $G$-orbit. The point is to find a small $G$ because in this case the group action can explain an otherwise mysteriously small denominator in the average of~$f$ across all of $\mathcal{S}$. If $f(\mathcal{S}) \subseteq \mathbb{Z}$ then of course this denominator is bounded by the order of $G$.

So to give a better proof of Theorem~\ref{thm:main} we could construct an action of~$C_3$, the cyclic group of order $3$, on $\mathcal{OT}(\lambda,k+2n)$ such that $\mathrm{wt}(\cdot)$ is homomesic with respect to this action. We might want this action to be free; if so, we should assume $n \geq 2$ so that $\#\mathcal{OT}(\lambda,k+2n)$ is always divisible by $3$. If the action is free then the size of each orbit must be $3$, so we only need for each $C_3$-orbit $\mathcal{O}$ that $\sum_{T \in \mathcal{O}} \mathrm{wt}(T) =   \frac{1}{2}(4n^2 + 3k^2 +8 kn + 2n +3k)$. We can make this request for a free cyclic action more precise via the relationship between oscillating tableaux and perfect matchings. We need to establish some preliminaries before we can state this relationship in Theorem~\ref{thm:rs}.

A \emph{perfect matching of $[2n]$} is a partition of the set $[2n] := \{1,2,\ldots,2n\}$ into pairs. Denote the set of perfect matchings of $[2n]$ by~$\mathcal{M}_{n}$. Let $M \in \mathcal{M}_n$ and let~$p = \{\alpha,\beta\}, q = \{\gamma,\delta\} \in M$. Assume that~$\alpha<\beta$, $\gamma<\delta$ and~$\beta < \delta$. Then we say that $\{p,q\}$ is a \emph{crossing} of $M$ if~$\alpha < \gamma < \beta < \delta$. We say it is a \emph{nesting} if~$\gamma < \alpha < \beta <\delta$. Finally, we say it is an \emph{alignment} if it is neither a crossing nor a nesting. We will denote the number of crossings, nestings, and alignments of~$M$ by $\mathrm{cr}(M)$, $\mathrm{ne}(M)$, and $\mathrm{al}(M)$ respectively.

Recall that a \emph{Dyck path of semilength $n$} is a lattice path in $\mathbb{Z}^2$ from $(0,0)$ to $(n,n)$ consisting of $n$ up steps of the form $(1,0)$ and $n$ down steps of the form $(0,1)$ that never goes below the diagonal $x=y$. Let~$\mathcal{D}_n$ denote the set of Dyck paths of semilength $n$. We identify a Dyck path $D \in \mathcal{D}_n$ with its Dyck word, which is the binary word of length $2n$ whose $i$th letter is $1$ if the $i$th step of $D$ is up and $0$ if it is down.  The \emph{area} of $D \in \mathcal{D}_n$, denoted $\mathrm{area}(D)$, is the number of complete boxes of unit size below $D$ and above the diagonal $x = y$. So for instance $\mathrm{area}(101010) = 0$, $\mathrm{area}(101100) = 1$, and $\mathrm{area}(111000) = 3$.

For $M \in \mathcal{M}_n$ and $\{i,j\} \in M$, if $i < j$ then we say $i$ is an \emph{opener} in~$M$ and $j$ is a \emph{closer} in~$M$. There is a natural surjection $\mathcal{M}_n \twoheadrightarrow \mathcal{D}_n$ given by~$M \mapsto D_M$ where the $i$th letter of $D_M$ is~$1$ if $i$ is an opener in $M$ and~$0$ if~$i$ is a closer in $M$. There is also a natural surjection $\mathcal{OT}(\emptyset,2n) \twoheadrightarrow \mathcal{D}_n$ given by $T = (\lambda^{0},\ldots,\lambda^{2n}) \mapsto D_T$ where the $i$th letter of $D_T$ is~$1$ if $\lambda^{i-1} \subdot \lambda^{i}$ and $0$ if $\lambda^{i} \subdot \lambda^{i-1}$.

\begin{thm} \label{thm:rs}
There is a bijection $RS\colon \mathcal{M}_n \xrightarrow{\sim} \mathcal{OT}(\emptyset,2n)$ such that $D_M = D_{RS(M)}$  for all~$M \in \mathcal{M}_n$.
\end{thm}

This bijection is apparently due to Richard Stanley; but we call it $RS$ also because it is an extension of Robinson-Schensted insertion. A description of the bijection is given in~\cite[Lemma 8.3]{sundaram}. See also~\cite[\S5]{chen}. In light of Theorem~\ref{thm:rs} it makes sense to ask how the statistic weight for oscillating tableaux translates to perfect matchings. That question is answered by the following proposition.

\begin{prop} \label{prop:stats}
\hspace{2em}
\begin{enumerate}[(a)]
\item For all $M \in \mathcal{M}_n$ we have $\mathrm{al}(M) = \binom{n}{2} - \mathrm{area}(D_M)$.
\item For all $T \in \mathcal{OT}(\emptyset,2n)$ we have $\mathrm{wt}(T) = 2\cdot\mathrm{area}(D_T) + n$.
\end{enumerate}
Thus for all $M \in \mathcal{M}_n$ we have $\mathrm{wt}(RS(M)) = n + 2\left(\binom{n}{2}-\mathrm{al}(M)\right)$.
\end{prop}

\noindent \emph{Proof}: Let $D \in \mathcal{D}_n$. For $1 \leq i \leq 2n$  define $b_i(D)$ to be the number of $1$s minus the number of $0$s in the subword consisting of the first $i$ letters of $D$. For $1 \leq i \leq n$ define $a_i(D) := b_j$ where $j$ is the position of the $i$th $0$ in $D$. If  $T = (\lambda^{0},\ldots,\lambda^{2n})$ then~$b_i(D_T) = |\lambda^{i}|$ and therefore $\mathrm{wt}(T)  = \sum_{i=1}^{2n} b_i(D_T)$. We also claim that $\mathrm{cr}(M) + \mathrm{ne}(M) = \sum_{i=1}^{n} a_i(D_M)$. This is because if the position of the $i$th $0$ in $D_M$ is $\beta$ and $\{\alpha,\beta\} \in M$ then
\begin{align*}
a_i(D_M) &= \# \{\{\gamma,\delta\} \in M\colon \gamma < \delta \textrm{ and } \gamma < \beta\} - \#\{\{\gamma,\delta\} \in M\colon \gamma < \delta \leq \beta\} \\
&= \# \{\{\gamma,\delta\} \in M \colon \alpha < \gamma < \beta < \delta \textrm{ or } \gamma < \alpha < \beta < \delta\}.
\end{align*}
Since $\mathrm{cr}(M) + \mathrm{ne}(M) + \mathrm{al}(M) = \binom{n}{2}$ for any $M \in \mathcal{M}_n$, we get $\mathrm{al}(M) = \binom{n}{2} - \sum_{i=1}^{n} a_i(D_M)$. So the proposition amounts to showing $\sum_{i=1}^{n} a_i(D) = \mathrm{area}(D)$ and $\sum_{i=1}^{2n} b_i(D) = 2\cdot\mathrm{area}(D) + n$ for any $D \in \mathcal{D}_n$. But this is true because it is true when $D = 1010\cdots10$ and it remains true after replacing a consecutive subword of $01$ by $10$ in any Dyck word. $\square$

It is not hard to convince oneself that among all $M \in \mathcal{M}_n$ and all $\{p,q\} \subseteq M$ it is equally likely for $\{p,q\}$ to be a crossing, nesting, or alignment. Thus the average number of alignments among all matchings is given by~$\frac{1}{\#\mathcal{M}_n} \sum_{M \in \mathcal{M}_n} \mathrm{al}(M) = \frac{1}{3}\binom{n}{2}$. This fact together with Theorem~\ref{thm:rs} and Proposition~\ref{prop:stats} yields Theorem~\ref{thm:main} in the case where $\lambda = \emptyset$. It also suggests that to find a~$C_3$ action on~$\mathcal{OT}(\lambda,k+2n)$ that exhibits homomesy with $\mathrm{wt}(\cdot)$ we could start by looking for a $C_3$ action on $\mathcal{M}_n$ that exhibits homomesy with~$\mathrm{al}(\cdot)$.

Presumably it would be easy to find such an action if the statistics $\mathrm{cr}(\cdot)$, $\mathrm{ne}(\cdot)$, and $\mathrm{al}(\cdot)$ were symmetrically distributed. However, this is not the case. That is, it is not true in general that
\begin{align} \label{eqn:dist}
 \sum_{M \in \mathcal{M}_{n}} x_1^{\mathrm{cr}(M)}x_2^{\mathrm{ne}(M)}x_3^{\mathrm{al}(M)} = \sum_{M \in \mathcal{M}_{n}} x_{\omega(1)}^{\mathrm{cr}(M)}x_{\omega(2)}^{\mathrm{ne}(M)}x_{\omega(3)}^{\mathrm{al}(M)}
 \end{align}
for all $\omega \in \mathfrak{S}_3$, the symmetric group on $3$ letters. Nevertheless, de M\'{e}dicis and Viennot~\cite{demedicis} define a certain involution, let us call it $\sigma\colon \mathcal{M}_n \to \mathcal{M}_n$, such that $\mathrm{cr}(M) = \mathrm{ne}(\sigma(M))$ and~$\mathrm{ne}(M) = \mathrm{cr}(\sigma(M))$ for all $M \in \mathcal{M}_n$.\footnote{There is another natural involution $c\colon \mathcal{OT}(\emptyset,2n) \to \mathcal{OT}(\emptyset,2n)$ that can be transferred to perfect matchings via $RS$ given by $c(\lambda^{0},\ldots,\lambda^{2n}) := ((\lambda^{0})',(\lambda^{1})',\ldots,(\lambda^{2n})')$. It is this conjugation symmetry that Chen~et.~al.~\cite{chen} exploit in studying maximal crossing and maximal nesting numbers in perfect matchings and set partitions. Note that $c$ is not the same as $\sigma$. One way to understand their relationship is by restricting to permutations. The set of $M \in \mathcal{M}_n$ with $D_M = 11\cdots 100\cdots 0$ is in bijection with $\mathfrak{S}_n$ by the map $M \mapsto \omega_M$ where $\omega_M(j) = i$ if $\{i,j+n\} \in M$. For all such $M$ we have $\omega_{c(M)} = \omega_{M}^{-1}$ while $\omega_{\sigma(M)} = \omega_{M}^{\mathrm{rev}}$. At least since the fundamental work of Simion and Schmidt~\cite{simion} it has been known that inversion and reversing are symmetries of the permutation pattern containment poset, and from the work of Smith~\cite{smith} it follows that these involutions generate the group of all symmetries of this poset (which is isomorphic to the dihedral group of order eight). Researchers such as Bloom and Elizalde~\cite{bloom2} have considered pattern containment among perfect matchings. It would be interesting to describe all the symmetries of the perfect matching pattern containment poset, especially in relation to the maps $c$ and $\sigma$.} (See~\cite{kasraoui} for a description of this involution in English as well as a generalization to arbitrary set partitions.) This shows that~(\ref{eqn:dist}) holds whenever $\omega(3)=3$. It also suggests that $\mathrm{cr}(\cdot)$, $\mathrm{ne}(\cdot)$, and $\mathrm{al}(\cdot)$ may posses a $\mathfrak{S}_3$ symmetry in spite of not being symmetrically distributed. This is because we could find some $\tau\colon \mathcal{M}_n \to \mathcal{M}_n$ of order $3$ such that $\mathrm{al}(\cdot)$ is homomesic with respect to~$\tau$. Then if we set~$G := \langle \sigma, \tau \rangle$ we would get that all three of the statistics $\mathrm{cr}(\cdot)$, $\mathrm{ne}(\cdot)$, and~$\mathrm{al}(\cdot)$ are homomesic with respect to the action of $G$. And if $(\sigma\tau)^2 = \mathrm{id}$ then $G \simeq \mathfrak{S}_3$. To our knowledge, such a $\mathfrak{S}_3$ action would be the first instance of homomesy involving a non-cyclic group. Constructing such a~$\tau$ remains an open problem. But it is useful to know, as Theorem~\ref{thm:main} suggests, that~$\tau$ ought to make sense at the more general level of oscillating tableaux of arbitrary shape as well.

\bibliography{osc}{}
\bibliographystyle{plain}

\end{document}